\documentclass[12pt]{article}

\usepackage{amssymb}

\def\M{{\cal M}}
\def\oM{{\overline{\cal M}}}
\def\Q{{\mathbb Q}}
\def\A{{\cal A}}
\def\qed{{\hfill $\diamondsuit$}}
\def\cL{{\cal L}}
\def\CP{{{\mathbb C}{\rm P}}}
\def\Aut{{\rm Aut}}

\usepackage{theorem}

\newtheorem{theorem}{Theorem}
\newtheorem{proposition}{Proposition}[section]
\newtheorem{corollary}[proposition]{Corollary}

\newtheorem{conjecture}[proposition]{Conjecture}
{\theorembodyfont{\rmfamily}
\newtheorem{definition}[proposition]{Definition}
\newtheorem{example}[proposition]{Example}
\newtheorem{remark}[proposition]{Remark}
\newtheorem{notation}[proposition]{Notation}
}

\include{epsf}

\title{An algebra of power series arising in
the intersection theory of moduli spaces of curves
and in the enumeration of ramified coverings of the sphere}

\author{Dimitri Zvonkine\thanks{
Institut f\"ur Mathematik,
Universit\"at Z\"urich,
Winterthurerstra\ss{}e 190,
CH-8057 Z\"urich.
E-mail: zvonkine@math.unizh.ch\,.\/\/
The author was partially suported by
EAGER - European Algebraic Geometry Research Training Network, 
contract No. HPRN-CT-2000-00099 (BBW) and by the Russian Foundation
of Basic Research grant 02-01-22004.}}

\date{\today}

\begin{document}

\maketitle

\begin{abstract}
A bracket is a function that assigns a number to
each monomial in variables $\tau_0, \tau_1, \dots$.
We show that any bracket satisfying the string and
the dilaton relations gives rise to a power series
lying in the algebra $\A$ generated by the series 
$\sum n^{n-1} q^n/n!$ and $\sum n^n q^n /n!\,$. 

As a consequence, various series from $\A$
appear in the intersection theory of moduli spaces
of curves. 

A connection between the counting of
ramified coverings of the sphere and the intersection
theory on moduli spaces allows us to prove that
some natural generating functions enumerating the
ramified coverings lie, yet again, in $\A$.
As an application, one can find the asymptotic
of the number of such coverings as the number of
sheets tends to $\infty$.

We believe that the leading terms of the asymptotics
like that correspond to observables in 2-dimensional
gravity.
\end{abstract}

\tableofcontents

\section{Introduction}

Denote by $\A$ the subalgebra of the algebra of
power series in one variable, generated by the series
$$
\sum_{n \geq 1} \frac{n^{n-1}}{n!} \, q^n \quad \mbox{and} \quad
\sum_{n \geq 1} \frac{n^n}{n!} \, q^n.
$$
We wish to show that this algebra plays an important role
in the intersection theory of moduli spaces $\oM_{g,n}$
of stable curves and in the problem of enumeration of
ramified coverings of the sphere.

\paragraph{SECTION~\ref{Sec:A}} 
contains a more explicit description of the algebra $\A$. 
In this section we also prove some relations between $\A$
and the combinatorics of Cayley trees ($=$ trees with
numbered vertices).

\paragraph{SECTION~\ref{Sec:M_{g,n}}} 
is devoted to the intersection theory on moduli spaces.

\begin{notation} \label{Not:MandL}
We denote by $\M_{g,n}$ the moduli space of smooth genus $g$
curves with $n$ marked and numbered distinct points. 

Further, $\oM_{g,n}$ is the Deligne-Mumford compactification of
this moduli space; in other words, $\oM_{g,n}$ is the
space of {\em stable} genus $g$ curves with $n$ marked
points.

We also use the standard notation $\cL_i$ for the following
line bundle  over $\oM_{g,n}$: consider a point $x \in \oM_{g,n}$
and the corresponding stable curve $C_x$; then the fiber
of $\cL_i$ over $x$ is the cotangent line to $C_x$
at the $i$th marked point. These line bundles are called
{\em tautological}.
\end{notation}

We will need the first Chern
classes $c_1(\cL_i)$ of the above line bundles, more precisely
the expression
$$
\frac{1}{1-c_1(\cL_i)} = 1 + c_1(\cL_i) + c_1(\cL_i)^2 +
\dots \in H^*(\oM_{g,n}, \Q).
$$

Let $\beta \in H^*(\oM_{g,0},\Q)$
be any cohomology class of $\oM_{g,0}$. For any $n \geq 0$,
there is a forgetful map from $\oM_{g,n}$ onto $\oM_{g,0}$,
forgetting the marked points and
contracting the components of the curve that have become
unstable. By abuse of notation, 
the pull-back of $\beta$ to $\oM_{g,n}$ by this
map will be again denoted by $\beta$.

\begin{theorem}\label{Thm:integral} 
For any $g \geq 2$ and $\beta \in H^*(\oM_{g,0}, \Q)$,
the power series
$$
F_{g,\beta} (q) =
\sum_{n \geq 0} \frac{q^n}{n!}
\int_{\oM_{g,n}} \frac{\beta}{(1-c_1(\cL_1)) \dots (1-c_1(\cL_n))}
$$
lies in the algebra $\A$.
\end{theorem}

In Section~\ref{Ssec:exceptions} we describe what happens
for genus~0 and~1, when the moduli space $\oM_{g,0}$
does not exist.

Later we will give an important generalization of this theorem
(Theorem~\ref{Thm:gen}).

To prove Theorem~\ref{Thm:integral} 
we first establish the analogs of the
``string equation'' and the ``dilaton equation'' for the
integrals involved in the right-hand part. Then we interpret
these equations as operations on graphs with numbered
vertices and use combinatorial results on these graphs to
prove the theorem. It is important to note that the 
second part of the proof relies {\em only} on the
string and dilaton relations.

\paragraph{SECTION~\ref{Sec:coverings}}
is devoted to the problem of enumerating the ramified
coverings of the sphere with specified ramification types.

Consider a nonconstant holomorphic map $f : C \rightarrow \CP^1$
of degree $n$ from a smooth complex curve $C$ to the Riemann sphere. 
Such maps will be called {\em ramified coverings} with $n$
sheets.

A {\em ramification point} of $f$ is a point of the target Riemann
sphere that has less than $n$ distinct preimages. 

For each ramification point $y$ of a ramified covering, we are
going to single out several {\em simple} preimages of $y$.

\begin{definition} \label{Def:markedcov}
A {\em marked ramified covering} is a ramified covering
with a choice, for every ramification point $y$, of a subset
of the set of simple preimages of $y$.
\end{definition}

Consider a partition $\mu = 1^{a_1} 2^{a_2} \dots$ of an integer 
$m \leq n$. Here we use multiplicative notation for
partitions: the partition $\mu$ contains $a_1$ parts equal to
$1$, $a_2$ parts equal to $2$, and so on, $\sum i a_i = m$. 

Suppose that a point $y \in \CP^1$ has $a_1$ 
marked simple preimages, $a_2$ double preimages, and so on. 
(Consequently, $y$ also has $n-m$ unmarked simple preimages.)
We then say that $y$ is a ramification point of $f$ of {\em multiplicity}
$r=a_2+2a_3 + 3a_4+\dots$ and of {\em ramification type}
$\mu = 1^{a_1}2^{a_2} \dots$. Sometimes the number $r$ will
also be called the {\em degeneracy} of the partition $\mu$.

\begin{definition}
\label{Def:Hurwitz}
A {\em Hurwitz number} $h_{n; \mu_1, \dots, \mu_k}$ 
is the number of connected $n$-sheeted marked ramified coverings of $\CP^1$
with $k$ ramification points, whose ramification types
are $\mu_1, \dots, \mu_k$. Every such covering is counted with weight
$1/|\Aut|$, where $|\Aut|$ is the number of automorphisms of the
covering.
\end{definition}

Note that the genus $g$ of the covering surface can be reconstituted
from the data $n; \mu_1, \dots, \mu_k$ using the Riemann-Hurwitz
formula: if the degeneracy of $\mu_i$ equals $r_i$, then
$$
2-2g = 2n-r_1- \dots - r_k.
$$

\bigskip

Fix $k$ nonempty partitions $\mu_1, \dots, \mu_k$ with
degeneracies $r_1, \dots, r_k$. Let $r$
be the sum $r = r_1 + \dots + r_k$. 

\begin{notation}
\label{Not:Hurwitz}
Denote by $h_{g,n;\mu_1, \dots, \mu_k}$ the number of $n$-sheeted
marked ramified coverings of $\CP^1$ by a genus $g$ surface,
with $k$ ramification points of types $\mu_1, \dots, \mu_k$
and, in addition, $2n+2g-2-r$ simple ($=$ of multiplicity~1)
ramification points. Each covering is counted with weight
$1/|\Aut|$.
\end{notation}

\begin{theorem} \label{Thm:Hurwitz}
Fix any $g \geq 0$, $k \geq 0$. If $g=1$, we suppose that
$k \geq 1$. Then for any partitions
$\mu_1, \dots, \mu_k$, the series
$$
H_{g;\mu_1, \dots, \mu_k} (q) =
\sum_{n \geq 1} \frac{h_{g,n;\mu_1, \dots, \mu_k}}{n!} 
\, q^n
$$
lies in the algebra $\A$.
\end{theorem}

This theorem is proved by induction on the number $k$ of
nonsimple ramification points. Curiously, the base of induction
(the case $k=1$) is the most difficult part of the proof. We 
prove it using a generalization of Theorem~\ref{Thm:integral}
and the Ekedahl-Lando-Shapiro-Vainshtein (or ``ELSV'') 
formula~\cite{ELSV} relating Hurwitz numbers to 
integrals over moduli spaces of stable curves.

Among other things, this theorem allows one to find the
asymptotic of the coefficients of $H_{g;\mu_1, \dots, \mu_k}$ as 
$n \rightarrow \infty$, knowing only the several first
coefficients.

\paragraph{SECTION~\ref{Sec:conclusion}} describes a
relation between the asymptotic of Hurwitz numbers 
and 2-dimensional gravity. 

We also compare the enumerative problems concerning ramified
coverings of the sphere and of the torus. While the former
are related to the intersection theory on $\oM_{g,n}$
and give rise to the algebra $\A$, the latter are related
with volumes of spaces of abelian differentials on Riemann
surfaces and give rise to the algebra of quasi-modular forms.

\paragraph{Acknowledgments}

The author is grateful to
J.-M.~Bismut, F.~Labourie, M.~Kontsevich, S.~Lando, S.~Natanzon, Ch.~Okonek,
A.~Okounkov, D.~Panov, J.-Y. Welschinger,
D.~Zagier, A.~Zorich and A.~Zvonkin for useful
discussions. I also thank for their interest 
the participants of the mathematical physics seminar at the ETH Z\"urich
and of the mathematical seminar at the ENS Lyon, as
well as the participants of the Luminy conference on
billiards ans Teichm\"uller spaces.

This work was partially supported by
EAGER - European Algebraic Geometry Research Training Network, 
contract No. HPRN-CT-2000-00099 (BBW) and by the RFBR grant
02-01-22004.

\paragraph{List of notations.} Here we summarize some
notation that we use consistently throughout the
paper. We have tried to avoid using the same letter
in different contexts unless it emphasizes a relation
between the two contexts that appears in some of
the proofs.

\bigskip

\noindent
\begin{tabular}{lp{29em}}
$n$ & The number of vertices in a graph. The number
of sheets of a covering.
The power of the variable $q$ in generating series. 
The number of marked points on a Riemann surface is
sometimes $n$ and sometimes $n-r$.\\
$q$ & The variable in generating series (to a sequence
$s_n$ we usually assign the series $\sum s_n q^n/n!$).\\
$g$ & The genus of a Riemann surface.\\
$\mu$ & A partition.\\
$p$ & The number of parts of a partition $\mu$.\\
$a_i$ & The number of parts of a partition $\mu$ that are
equal to $i$.\\
$b_i$ & The parts of a partition $\mu$
are denoted by $b_1, \dots, b_p$.\\
$r$ & The degeneracy of a partition $\mu$ defined
by $r= \sum (b_i-1)$. The multiplicity
of a ramification point.\\
$k$ & The number of partitions. If $k > 1$, the
partitions are denoted by $\mu_1, \dots, \mu_k$,
their degeneracies by $r_1, \dots, r_k$, while
$r$ is the total degeneracy $r = \sum r_i$.\\
$c(n)$ & The number of simple ramification points
in a ramified covering.\\
$\beta$ & A cohomology class on the moduli space
$\oM_{g,p}$~.\\
$b$ & We often take $\beta$ to be a 
cohomology class of pure degree $2b$.\\
$\psi_i$ & The first Chern class $c_1(\cL_i)$
of the line bundle $\cL_i$.\\
$d_i$ & The power of the class $\psi_i$ in the
intersection numbers we consider. We also consider
graphs whose $i$th vertex has valency $d_i+1$
for each~$i$.
\end{tabular}

\section{The algebra $\A$ of power series}
\label{Sec:A}

The algebra of power series
$$
\A = \Q \left[ \sum_{n \geq 1} \frac{n^{n-1}}{n!} \, q^n,
\sum_{n \geq 1} \frac{n^n}{n!} \, q^n \right]
$$
plays a central role in this paper.
Here we give an explicit description of $\A$ and show its
relation with the combinatorics of Cayley trees. Many of
the results below are known, but have probably never
been put together. As far as we know, the algebra $\A$
itself was first discovered by D.~Zagier several years ago
(unpublished), and then independently 
introduced in our paper~\cite{Zvonkine},
where most of the results of Section~\ref{Ssec:compA} are
given.

\subsection{How to make computations in $\A$}
\label{Ssec:compA}

Denote by $Y$ and $Z$ the generators of $\A$
\begin{eqnarray*}
Y &=& \sum_{n \geq 1} \frac{n^{n-1}}{n!} \, q^n, \\
Z &=& \sum_{n \geq 1} \frac{n^n}{n!} \, q^n.
\end{eqnarray*}
Denote by $D$ the differential operator 
$D = q \, \frac{\partial}{\partial q}$. Thus $Z = DY$.

Note that both $Y$ and $Z$ have a radius of convergence
of $1/e$. Therefore the same is true of all series
in $\A$. The function $Y(q)$, more precisely, $-Y(-q)$,
was considered by J.~H.~Lambert~\cite{Lambert} in
1758 (we thank N. A'Campo for this reference).

\begin{proposition}\label{Prop:Y=qe^Y}
We have
$$
Y = q e^Y.
$$
\end{proposition}

\paragraph{Proof.} $Y$ is the exponential generating series
for rooted Cayley trees (Definition~\ref{Def:Cayley}). 
Therefore $e^Y$ is the exponential generating series for 
forests of rooted Cayley trees. Add
a new vertex $*$ to such a forest and join
$*$ to the root of each tree.
We obtain a Cayley tree with root $*$. This
operation is a one-to-one correspondence, hence 
$Y = q e^Y$. \qed

\begin{corollary}
On the disc $|q| < 1/e$, the function $Y(q)$ is the inverse of
the function $q(Y) = Y/e^Y$.
\end{corollary}

\begin{proposition} \label{Prop:(1-Y)(1+Z)}
We have $(1-Y)(1+Z) = 1$.
\end{proposition}

\paragraph{Proof.} 
$$
Z = DY = D(q e^Y) = q e^Y + q e^Y DY = q e^Y (1+Z)
= Y (1+Z).
$$
Hence $(1-Y)(1+Z) = 1$. \qed

\begin{corollary}
As an abstract algebra, $\A$ is isomorphic to
$\Q [X, X^{-1}]$, where $X=1-Y$.
\end{corollary}

Now we can show how to make computations in $\A$. For instance,
let us prove that 
$$
Y^2 = \sum_{n \geq 1} \frac{2(n-1)n^{n-2}}{n!} \, q^n
= 2 \sum_{n \geq 2} \frac{n^{n-3}}{(n-2)!} \, q^n.
$$
Indeed, $D(Y^2) = 2Y DY = 2YZ = 2(Z-Y)$, which determines
$Y^2$ up to its free term, but the free term obviously vanishes.

More generally, the powers of $Y$ are given by the following 
proposition.

\begin{proposition} \label{Prop:Y^k} We have
$$
Y^k = \sum_{n \geq 1} \frac{k(n-1) \dots (n-k+1) n^{n-k}}{n!} \, q^n
= k \sum_{n \geq k} \frac{n^{n-k-1}}{(n-k)!} \, q^n.
$$
\end{proposition}

\paragraph{Proof.} 
Induction on $k$. For $k=1$ the assertion is true. 
To go from $k$ to $k+1$, one checks that the equality
$$
D \left(
\frac{Y^{k+1}}{k+1} - \frac{Y^k}k
\right) = (Y^k - Y^{k-1}) DY = Y^{k-1} (Y-1)Z = -Y^k 
$$
is compatible with our expressions for $Y^k$ and $Y^{k+1}$.
The constant term of $Y^{k+1}$ vanishes, thus $Y^{k+1}$
is uniquely determined by $D(Y^{k+1})$. \qed

\bigskip

Now we study the powers of $Z$.

\begin{definition} \label{Def:A_n}
Denote by $A_n$ the sequence of integers
$$
A_n = \sum_{\stackrel{\scriptstyle p+q=n}{p,q \geq 1}}
\frac{n!}{p! \, q!} \, p^p \, q^q.
$$
\end{definition}
Its first terms are $0, 2, 24, 312, 4720, \dots$.
We have
$$
Z^2 = \sum_{n \geq 1} \frac{A_n}{n!} \, q^n.
$$
One can show that 
$$
A_n = n! \, \sum_{k=0}^{n-2} \frac{n^k}{k!}   
\sim \sqrt{\pi/2} \; n^{n+\frac12}.
$$
As far as we know, there is no simple expression for
the powers of $Z$. However, we can prove that they are
linear combinations of the series
$$
D^k Z = \sum_{n \geq 1} \frac{n^{n+k}}{n!} \, q^n
$$
and
$$
D^k (Z^2) = \sum_{n \geq 1} \frac{n^k A_n}{n!} \, q^n \, .
$$

\begin{proposition} \label{Prop:Z^k}
For any integer $k \geq 0$, the power series $D^k Z$ and $D^k (Z^2)$
are polynomials in $Z$ with positive integer coefficients,
of degrees  $2k+1$ and $2k+2$ respectively.
\end{proposition}

\paragraph{Proof.}
Applying $D$ to both sides of the equality $(1-Y)(1+Z)=1$
we get
$$
-Z(1+Z) + (1-Y) \cdot DZ = 0.
$$
Thus
$$
DZ = \frac{Z(1+Z)}{1-Y} = Z(1+Z)^2.
$$
Hence
$$
D(Z^2) = 2Z^2(1+Z)^2.
$$
Now we proceed by induction on $k$. \qed

\begin{corollary} \label{Cor:Z^k}
For any positive integer $k$, the power series $Z^k$ is a
linear combination with integer coefficients
of the first $k$ series from the list
$Z, Z^2, DZ, D(Z^2), D^2 Z, D^2 (Z^2), \dots$.
\end{corollary}

From Proposition~\ref{Prop:Y^k} and Corollary~\ref{Cor:Z^k}
we deduce the following theorem.

\begin{theorem} {\rm \cite{Zvonkine}}
The algebra $\A$ is spanned over $\Q$ by the power series
$$
1, 
\qquad \quad 
\sum_{n \geq 1} \frac{n^{n+k}}{n!} q^n,
\quad k \in {\mathbb Z}, 
\qquad \quad 
\sum_{n \geq 1} \frac{n^k A_n}{n!} q^n,
\quad k \in {\mathbb N}.
$$
\end{theorem}

Note that the Stirling formula together with the asymptotic
for the sequence $A_n$ allows one to determine the
leading term of the asymptotic for the coefficients of any 
series in $\A$. We have
$$
\frac{n^n}{n!} \sim \frac1{\sqrt{2\pi n}} \; e^n, 
\qquad
\frac{A_n}{n!} \sim \frac12 \, e^n.
$$

Note also that if, for some series $F \in \A$,
we know in advance its degree in $Y$ and in $Z$,
then we can reconstitute the series $F$ using only
a finite number of its initial terms. 

Combining both remarks, we see that initial terms of the sequence 
of coefficients of $F$ determine the asymptotic of the sequence!

\subsection{Dendrology}
\label{Ssec:dendrology}

\begin{definition} \label{Def:Cayley}
A {\em Cayley tree} is a tree with numbered vertices.
\end{definition}

It is well-known (Cayley theorem) that there are $n^{n-2}$ 
Cayley trees with $n$ vertices. Note that the corresponding
exponential generating function
$$
\sum_{n \geq 1} \frac{n^{n-2}}{n!} \, q^n
$$
lies in the algebra $\A$.

Consider a Cayley tree $T$ with two marked vertices $a$ and $b$. 
Denote by $l(T)$ the distance between these vertices, 
i.e., the number of edges in the shortest path joining them.

\begin{definition}
Denote by $m_{n,k}$ and $p_{n,k}$ the sums 
$$
m_{n,k} = \sum_{T} l(T)^k, \qquad 
p_{n,k} = \sum_{T} \frac{l(T) (l(T)-1) \dots (l(T)-k+1)}{k!}
$$
where the sum is taken over all Cayley trees $T$ with $n$
vertices, two of which are marked.
\end{definition}

For instance, $m_{2,1} = p_{2,1} = 2$. Note that if we consider $l(T)$ as
a random variable, then $m_{n,k}$ is its $k$th moment.

\begin{theorem} \label{Thm:dendrology}
For any $k$, the power series
$$
\sum_{n \geq 1} \frac{m_{n,k}}{n!} \, q^n
\quad \mbox{and} \quad
\sum_{n \geq 1} \frac{p_{n,k}}{n!} \, q^n
$$
lie in $\A$.
\end{theorem}

\begin{example}
It follows from the proof below 
that $p_{n,1} = m_{n,1} = A_n$. This number
is called {\em the total height of Cayley trees} and
was introduced in~\cite{RioSlo}.
\end{example}

\paragraph{Proof of Theorem~\ref{Thm:dendrology}.}
It is sufficient to prove the theorem for $p_{n,k}$.

Fix $k$. There is a natural bijection between the following
sets of objects. 

$E_n$ is the set of Cayley trees with
$n$ vertices, on which one has labeled two vertices
by $a$ and $b$ and chosen $k$ distinct edges on the
shortest path from $a$ to $b$. The number of elements
of $E_n$ equals $p_{n,k}$.

$F_n$ is the set of ordered $(k+1)$-tuples of trees with $n$ vertices
in whole; the vertices are numbered from $1$ to $n$ and,
in addition, two vertices $a_i$ and $b_i$, $1 \leq i \leq k+1$,
are marked on each tree.

The bijection is established in the following way: taking
a forest from the set $F_n$ we draw new edges joining the
vertices $b_1$ to $a_2$, then $b_2$ to $a_3$, \dots,
and finally $b_k$ to $a_{k+1}$. We obtain a tree with $k$
marked edges lying on the path between $a_1$ and $b_{k+1}$,
i.e., a tree from the set $E_n$.

Now, the trees with two marked vertices are enumerated by
the series $Z$, therefore the exponential generating series
for the sequence $|F_n|$ is $Z^{k+1}$.
\qed

\section{Intersection theory on $\oM_{g,n}$}
\label{Sec:M_{g,n}}

\subsection{The bracket 
$\left< \beta \, \tau_{d_1} \dots \tau_{d_n} \right>$}
\label{Ssec:bracket}

Recall that $\oM_{g,n}$ is the moduli space of stable curves
and $\cL_i$ the tautological line bundles (Notation~\ref{Not:MandL}).
Further, $\psi_i = c_1(\cL_i)$.
Let $g \geq 2$ and denote by $\beta \in H^*(\oM_{g,0}, \Q)$ a 
cohomology class of $\oM_{g,0}$ as well as its pull-backs to the 
spaces $\oM_{g,n}$, $n \geq 0$, under the projections forgetting the
marked points. 

\begin{notation}
For $n \geq 0$, 
we denote by $\left< \beta \, \tau_{d_1} \dots \tau_{d_n} \right>$
the integral
$$
\left< \beta \, \tau_{d_1} \dots \tau_{d_n} \right>
=
\int_{\oM_{g,n}} \beta \; \; \psi_1^{d_1} \dots  \psi_n^{d_n}.
$$
\end{notation}

\begin{remark}
If the cohomology class $\beta$ has pure degree $2b$, then
the bracket vanishes unless $b + \sum d_i = 3g-3+n$. If 
$\beta$ is odd, the bracket always vanishes.

In the case $\beta = 1$, our notation is compatible with
the standard notation for intersection numbers of the first Chern
classes of the bundles $\cL_i$ (see~\cite{Witten}).
\end{remark}

The power series from Theorem~\ref{Thm:integral},
which is one of the main objects of our study, can be
written in the form
$$
F_{g,\beta} (q) \; = \;
\sum_{n \geq 0}
\frac{q^n}{n!} 
\int_{\oM_{g,n}}
\frac{\beta}{(1-\psi_1) \dots (1-\psi_n)}
\; = \; 
\sum_{n \geq 0}
\frac{q^n}{n!} \sum_{d_1, \dots, d_n}
\left< \beta \; \tau_{d_1} \dots \tau_{d_n} \right>
 .
$$

\subsection{Simplest cases: $g=0$ and $g=1$}
\label{Ssec:exceptions}

Most of the theorems concerning the bracket
$\left< \beta \; \dots \; \right>$ have three
exceptional cases: $g=0$, $\deg \beta=0$
and $g=1$, $\deg \beta = 0 \mbox{ or } 2$.
We will describe these cases here, which will allow us
not to bother with them later and, at the same time, to give
the simplest examples.

Strictly speaking, these examples
are not totally compatible with our definitions,
because there is no space $\oM_{g,0}$ for $g=0,1$. However, 
in some sense, the cohomology classes below could be
seen as pull-backs from these non-existent spaces.
At least, the the bracket itself makes perfect sense. 

\paragraph{The case $g=0$, $\beta=1$.}
In this case, we have
$$
\left< \tau_{d_1} \dots \tau_{d_n} \right>_{g=0}
= 
\frac{(n-3)!}{d_1! \dots d_n!}
$$
if $\sum d_i = n-3$; otherwise the bracket vanishes.
Using the explicit values of the bracket
one gets
$$
F_{g=0,\beta=1} (q) = 
\sum_{n \geq 3} \frac{q^n}{n!}
\sum_{d_1, \dots, d_n}
\left< \tau_{d_1} \dots \tau_{d_n} \right>_{g=0}
=
\sum_{n \geq 3} \frac{n^{n-3}}{n!} \, q^n .
$$
The most logical thing to do is to add to this series
two terms, $q$ and $q^2/4$, that would correspond to
two non-existent moduli spaces $\oM_{0,1}$ and $\oM_{0,2}$.
Thus we obtain the series
$$
\sum_{n \geq 1} \frac{n^{n-3}}{n!} \, q^n \; \in \; \A \,.
$$

\paragraph{The case $g=1$, $\deg \beta = 0$.}
In this case, we have the following formula for
the values of the bracket. Denote by
$\sigma_j$ the $j$th elementary
symmetric function of $d_1, \dots, d_n$.
In other words, 
$$
(d-d_1) \dots (d-d_n) = d^n - \sigma_1 d^{n-1} +
\dots + (-1)^n \sigma_n.
$$
Then we have
$$
\left< \tau_{d_1} \dots \tau_{d_n} \right>_{g=1}
= 
\frac{1}{24 \, d_1! \dots d_n!}
\left(
n! - \sum_{j=2}^n (j-2)! (n-j)! \sigma_j
\right)
$$
if $\sum d_i = n$; otherwise the bracket vanishes.
From this formula (after a certain amount of calculations),
one can deduce
$$
\sum_{d_1, \dots, d_n}
\left< \tau_{d_1} \dots \tau_{d_n} \right>_{g=1}
=  \frac{1}{24} \left( \frac{A_n}{n} + n^{n-1} \right) \, ,
$$
with the sequence $A_n$ from Definition~\ref{Def:A_n}.
Thus
$$
F_{g=1,\beta=1} (q) = 
\sum_{n \geq 1} \frac{q^n}{n!}
\sum_{d_1, \dots, d_n}
\left< \tau_{d_1} \dots \tau_{d_n} \right>_{g=1}
=
\frac{1}{24} \;  
\sum_{n \geq 1} \left( \frac{A_n}{n} + n^{n-1} \right)
\frac{q^n}{n!} \; .
$$
The series $F_{g=1, \beta=1}$ does not lie in the algebra $\A$,
but the series 
$$
D F_{g=1,\beta=1}(q) = \frac{1}{24} \;  
\sum_{n \geq 1} \frac{A_n + n^n}{n!} \, q^n
$$ 
does. A generalization of
this fact is given in Theorem~\ref{Thm:gen}.

\paragraph{The case $g=1$, $\deg \beta = 2$.}
In this paragraph we denote by $\beta$ the unique 
2-cohomology class
of $\oM_{1,1}$ whose integral over the fundamental class
of $\oM_{1,1}$ is equal to 1. The other 2-cohomology
classes are proportional to $\beta$, because
the vector space $H^2(\oM_{1,1}, \Q)$ is 1-dimensional.

We have
$$
\left< \beta \; \tau_{d_1} \dots \tau_{d_n} \right>
=
\frac{(n-1)!}{d_1! \dots d_n!}
$$
if $\sum d_i = n-1$; otherwise the bracket vanishes.
The corresponding generating series equals
$$
F_{g=1,\beta} (q) = 
\sum_{n \geq 1} \frac{q^n}{n!}
\sum_{d_1, \dots, d_n}
\left< \beta \; \tau_{d_1} \dots \tau_{d_n} \right>
= \sum_{n \geq 1} \frac{n^{n-1}}{n!} \, q^n = Y(q).
$$
This series, of course, lies in $\A$.

\subsection{String and dilaton relations}
\label{Ssec:string}

The following proposition is well-known in the
case $\beta=1$ (see, for example,~\cite{Witten}). 
The proof in the general case is literally the same as for 
$\beta=1$, but we still give it here for completeness.

In the right-hand side of the string relation
below, we set, by convention, a bracket containing $\tau_{-1}$
to be $0$.

\begin{proposition}
\label{Prop:SandD}
For $n \geq 0$, we have
$$
\left<\beta \;
\tau_{d_1} \dots \tau_{d_n} \; \tau_0
\right>
=
\sum_{i=1}^n
\left< \beta \;
\tau_{d_1} \dots \tau_{d_i-1} \dots \tau_{d_n} 
\right> \qquad \mbox{\rm (string relation)}
$$
and
$$
\left< \beta \;
 \tau_{d_1} \dots \tau_{d_n} \; \tau_1
\right>
=
(2g-2+n)
\left< \beta \;
\tau_{d_1} \dots \tau_{d_n} 
\right> \qquad \mbox{\rm (dilaton relation)}.
$$
\end{proposition}

\begin{figure}[h]
\begin{center}
\
\epsfbox{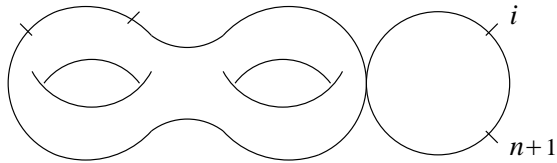}

\caption{A generic curve of the divisor $D_{i,n+1}$.}
\label{Fig:divisor}
\end{center}
\end{figure}

\paragraph{Proof.} Let $D_{i,n+1} \subset \oM_{g,n+1}$
be the divisor consisting of stable curves $C$ with the
following property: one of the irreducible 
components of $C$ is a sphere containing the marked points
$i$ and $n+1$ (and no other marked points) 
and exactly one node (Figure~\ref{Fig:divisor}).
We denote by $\Delta_i$ the $2$-cohomology class
Poincar\'e dual to $D_{i,n+1}$.

We will need to consider the line bundles $\cL_i$ both
on $\oM_{g,n}$ and on $\oM_{g,n+1}$. We momentarily
denote the former by $\cL'_i$, $1 \leq i \leq n$, the
latter retaining the notation $\cL_i$, $1 \leq i \leq n+1$.
Further, denote by $\psi'_i$
the first Chern class of the line bundle $\cL'_i$
on $\oM_{g,n}$ and, by abuse of notation, its pull-back
to $\oM_{g,n+1}$ by the map forgetting the $(n+1)$st marked
point. Let $\psi_i$ be the first Chern class of $\cL_i$
on $\oM_{g,n+1}$.

The forgetful map $\oM_{g,n+1} \rightarrow \oM_{g,n}$ identifies the line
bundles $\cL_i$ and $\cL'_i$ everywhere except
over the divisor $D_{i,n+1}$. Using this one can check that 
\begin{equation}
\psi_i = \psi'_i + \Delta_i. 
\end{equation}
Moreover, we have 
\begin{equation}
\psi_i \cdot \Delta_i = \psi_{n+1} \cdot \Delta_i = 0
\end{equation}
(because the line bundles $\cL_i$ and $\cL_{n+1}$
are trivial over $D_{i,n+1}$)
and
\begin{equation}
\Delta_i \cdot \Delta_j  = 0 \quad \mbox{for } i \not= j
\end{equation}
(because the divisors have an empty geometric intersection).

Let us first prove the string relation. We have
$$
\psi_i^d - (\psi'_i)^d \; \stackrel{(1)}{=} \;
\Delta_i \;
\left(
\psi_i^{d-1} + \dots + (\psi'_i)^{d-1}
\right) \;
\stackrel{(2)}{=} \;
\Delta_i \; (\psi'_i)^{d-1},
$$
where we set by convention $(\psi'_i)^{-1} = 0$.
Thus
$$
\psi_i^d \; = \;
(\psi'_i)^d + \Delta_i \; (\psi'_i)^{d-1}.
$$
It follows that
$$
\int_{\oM_{g,n+1}} 
\beta \; \psi_1^{d_1} \; \dots
\; \psi_n^{d_n} =
$$
$$
\int_{\oM_{g,n+1}} 
\beta \; 
\left[
(\psi'_1)^{d_1} + \Delta_1 \; (\psi'_1)^{d_1-1} 
\right]
\; \dots
\left[
(\psi'_n)^{d_n} + \Delta_n \; (\psi'_n)^{d_n-1} 
\right]
\;  \stackrel{(3)}{=}
$$
$$
\int\limits_{\oM_{g,n+1}} \!\!\!\!
\beta \; (\psi'_1)^{d_1} \; \dots
\; (\psi'_n)^{d_n} 
\; + \; 
\sum_{i=1}^n \; 
\int\limits_{\oM_{g,n+1}} \!\!\!\!
\beta \; (\psi'_1)^{d_1} \; \dots \;
\Delta_i \; (\psi'_i)^{d_i-1}
\; \dots
\; (\psi'_n)^{d_n} \; .
$$
The first integral is equal to $0$, because the integrand is
a pull-back from $\oM_{g,n}$. As for the integrals composing
the sum, we integrate the class $\Delta_i$ over the fibers of the
projection $\oM_{g,n+1} \rightarrow \oM_{g,n}$. This is equivalent
to restricting the integral to the divisor $D_{i,n+1}$,
which is naturally isomorphic to $\oM_{g,n}$. Finally, we obtain
$$
\int_{\oM_{g,n+1}} 
\beta \; \psi_1^{d_1} \; \dots
\; \psi_n^{d_n} =
\sum_{i=1}^n
\int_{\oM_{g,n}} 
\beta \; (\psi'_1)^{d_1} \; \dots \;
(\psi'_i)^{d_i-1}
\; \dots
\; (\psi'_n)^{d_n} \; .
$$

This proves the string relation.

Now let us prove the dilaton relation. We have
$$
\int_{\oM_{g,n+1}} 
\beta \; \psi_1^{d_1} \; \dots
\; \psi_n^{d_n} \; \psi_{n+1}
\stackrel{(1)}{=}
$$
$$
\int_{\oM_{g,n+1}}
\beta \; \biggl(\psi'_1 + \Delta_1 \biggr)^{d_1} \;
\dots \; \biggl(\psi'_n + \Delta_n \biggr)^{d_n} \; \psi_{n+1}
\stackrel{(2)}{=}
$$
$$
\int_{\oM_{g,n+1}}
\beta \; (\psi'_1)^{d_1} \; \dots
\; (\psi'_n)^{d_n} \; \psi_{n+1} =
$$
$$
(2g-2+n) \; \int_{\oM_{g,n}}
\beta \; (\psi'_1)^{d_1} \; \dots
\; (\psi'_n)^{d_n}.
$$
The last equality is obtained by integrating the factor 
$\psi_{n+1}$ over the fibers of the projection
$\oM_{g,n+1} \rightarrow \oM_{g,n}$.

This proves the dilaton relation.
\qed

\subsection{Graphs with numbered vertices}
\label{Ssec:graphs}

We are now going to interpret the string and the dilaton relations
as operations on graphs with numbered vertices. Each variable
$\tau_d$ will correspond to a vertex of valency $d+1$. The
string relation will allow us to erase a vertex of valency $1$,
simultaneously decreasing by $1$ the valency of its
neighboring vertex. The dilaton relation will allow us
to erase a vertex of valency $2$, merging together the
two edges that were adjacent to it. We will arrange
for our graphs to have $2g-2+n$ edges, which will
account for the factor $2g-2+n$ in the dilaton relation.
We will show that using such graphs one can define all possible
brackets $\left< \beta \; \tau_{d_1} \dots \tau_{d_n} \right>$
satisfying the string and dilaton relations.
Now let us give precise formulations.

\bigskip

The {\em valency} of a vertex in a graph is the number of
edges issuing from this vertex (a loop being counted twice).

We introduce two operations (S) and (D) on finite graphs. 

\begin{figure}[h]
\begin{center}
\
\epsfbox{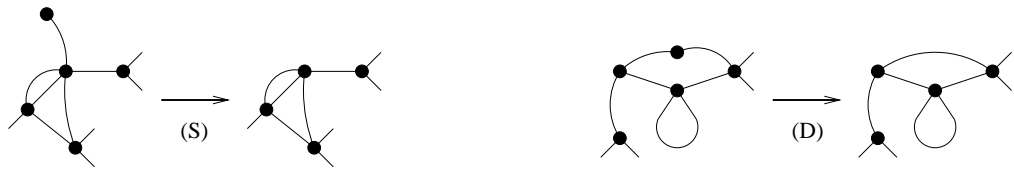}

\caption{Simplifying operations on graphs.}
\label{Fig:simplification}
\end{center}
\end{figure}

\begin{definition} \label{Def:SandD}
The {\em operation}~(S) consists in
erasing a vertex of valency $1$ and the edge adjacent to it.
The {\em operation}~(D) consists in erasing a vertex of valency $2$
and merging its two adjacent edges into one edge. These
operations are shown in Figure~\ref{Fig:simplification}.
\end{definition}

To define new brackets 
we will need the following type of graphs.

\begin{definition} \label{Def:decgraph}
We call an $n$-{\em decorated graph} a graph $G$ with $n+1$ vertices,
one of which is labeled by $*$ and the remaining ones are numbered
from $1$ to $n$. Moreover, we suppose that each 
connected component of $G$ either contains the vertex $*$
or has at least two independent cycles. Loops and
multiple edges are allowed.
\end{definition}

It is helpful to imagine that the factors $\beta$, $\tau_{d_1}$,
\dots $\tau_{d_n}$ from the bracket 
correspond, respectively, to the vertices
$*$, $1$, \dots, $n$. The vertex $*$ will soon be required
to have valency $g-1+\frac12 \deg \beta$. The vertex
number $i$, for each $i$, will soon be required to have
valency $d_i+1$.

Now we are going to apply the operations (S) and (D) to 
decorated graphs. The vertex $*$ is special: we will
erase neither $*$ itself nor the edges adjacent to it.

Let $G$ be an $n$-decorated graph. 

\begin{definition}
The {\em simplification of $G$} is the
graph $H$ obtained from $G$ by (i)~applying as many times as possible
the operations (S) and (D) without erasing the vertex $*$ or 
changing its valency; (ii)~forgetting the numbers of the vertices.

Conversely, the graph $G$ is called an {\em extension} of $H$.

A {\em simple} graph is any finite graph with
one vertex labeled by $*$ and several
other vertices, each of which has a valency $\geq 3$. 
\end{definition}

It is easy to see that the simplification of $G$ is well-defined 
and  unique, in other words, it does not depend
on the order in which one applies the operations (S) and (D).

Let $g \geq 2$ be an integer and $H$ a simple graph with
Euler characteristic 
$\chi(H) = |\mbox{vertices}| - |\mbox{edges}|$
equal to $3-2g$. For $n \geq 0$, we define a bracket
$\left< \tau_{d_1} \dots \tau_{d_n} \right>_H$ by the following
procedure.

\begin{definition} \label{Def:bracketH}
The value of the $H$-{\em bracket}
$\left< \tau_{d_1} \dots \tau_{d_n} \right>_H$ 
is equal to the weighted number of $n$-decorated graphs $G$ 
such that: (i)~the simplification of $G$ is $H$;
(ii)~for $1 \leq i \leq n$, the valency of the $i$th vertex of
$G$ equals $d_i+1$. The weight of a graph $G$ 
equals $1/\mbox{(number of its symmetries)}$.
\end{definition}

\begin{remark}
A {\em symmetry} of a decorated graph is a permutation of its half-edges
that preserves the vertex of each half-edge and does not split the
edges.
\end{remark}

\begin{proposition} \label{Prop:bracketH}
The bracket $\left< \; \cdot \; \right>_H$ satisfies the string and
dilaton relations of Proposition~\ref{Prop:SandD}.
\end{proposition}

\paragraph{Proof.} \ 

String relation. Suppose $d_{n+1} = 0$. Definition~\ref{Def:bracketH}
assigns to the list $d_1$, \dots, $d_n$, $0$ a set of
$(n+1)$-decorated graphs. Consider a graph $G$ from this set. 
Its $(n+1)$st vertex has valency $1$. First note that it 
cannot be joined to the vertex $*$. Indeed, otherwise the vertex 
number $n+1$ can never be erased by the operations (S) and (D), 
and therefore will remain as a vertex of $H$. But the
vertices of $H$ have valency at least $3$.
Thus the vertex number $n+1$ is joined to some vertex number $i$. 
Then, applying the operation (S) to the
$(n+1)$st vertex, we obtain an $n$-decorated graph $G'$
from the set assigned to the list $d_1, \dots,d_i-1, \dots, d_n$.
This operation is obviously a one-to-one correspondence between
the set of graphs $G$ corresponding to the list
$d_1, \dots, d_n,0$ and the disjoint union of the
sets of graphs $G'$ corresponding to the lists
$d_1, \dots, d_i-1, \dots, d_n$
for $1 \leq i \leq n$. This leads immediately to the string
relation for the bracket $\left< \; \cdot \; \right>_H$.

Dilaton relation. Suppose $d_{n+1}=1$. As above, consider an
$(n+1)$-decorated graph $G$ assigned to the list
$d_1, \dots, d_n, 1$. Its $(n+1)$st vertex has valency $2$.
Applying operation (D) to this vertex, we obtain a graph $G'$
from the set assigned to the list $d_1, \dots, d_n$.
This time, the operation is not a bijection, but becomes one
if we mark the edge of $G'$ on which the $(n+1)$st vertex of $G$
has disappeared. More precisely, we have a one-to-one
correspondence between the set of graphs $G$ assigned to the list
$d_1, \dots, d_n, 1$ and the set of pairs $(G',e)$,
where $G'$ is graph assigned to $d_1, \dots, d_n$ and
$e$ is its edge. Now, recall that the simple graph $H$,
and therefore the graphs $G$ and $G'$ as well,
have Euler characteristic $3-2g$. Thus each graph $G'$ has
$2g-2+n$ edges. This leads to the dilaton relation for
the bracket $\left< \; \cdot \; \right>_H$. \qed

\begin{remark}
The string and dilaton relations are {\em linear}, i.e., if several
different brackets satisfy these relations then so does
their arbitrary linear combination. Therefore 
Proposition~\ref{Prop:bracketH} immediately implies that 
the two relations are satisfied by any linear combination
$$
s_1 \left< \; \cdot \; \right>_{H_1} \;  + \; \dots \; 
+ \;  s_l \left< \; \cdot \; \right>_{H_l}
$$ 
for simple graphs $H_1, \dots, H_l$.
\end{remark}

Now we are going to prove that linear combinations of
brackets $\left< \; \cdot \; \right>_H$ for different
simple graphs $H$ represent {\em all possible} brackets
of nonnegative degree (see Definition~\ref{Def:degree} below)
that satisfy the string and the dilaton relations. In
particular, for any cohomology class $\beta$, the
bracket $\left< \beta \; \dots \; \right>$ can be
represented in that way.

\begin{figure}[h]
\begin{center}
\
\epsfbox{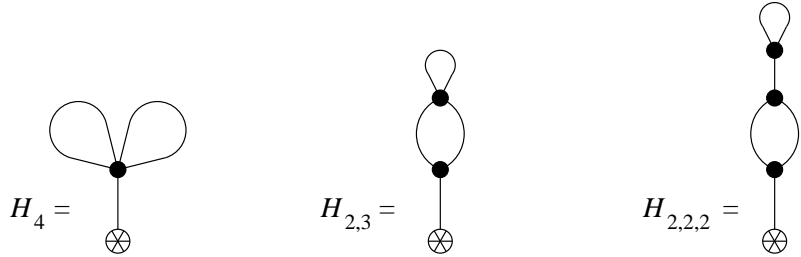}

\caption{The 3 graphs used to represent the bracket
$\left< \; \cdot \; \right>$ for $g=2$, $\beta=1$.}
\label{Fig:3graphs}
\end{center}
\end{figure}

\begin{example}
Consider the case $g=2$, $\beta = 1$. The corresponding
bracket is entirely determined by the string and dilaton relations
together with the ``initial conditions'' (we do not
explain here how these are found):
$$
\left< \tau_4 \right> = \frac1{1152}\, , \quad
\left< \tau_2 \tau_3 \right> = \frac{29}{5760} \, , \quad
\left< \tau_2 \tau_2 \tau_2 \right> = \frac7{240} \, .
$$
For each of the three brackets above, we choose a simple
graph that will represent it, for instance those shown in
Figure~\ref{Fig:3graphs}. Now we have
$$
\left< \; \cdot \; \right>_{g=2} 
=
8 \cdot \frac1{1152} \left< \; \cdot \; \right>_{H_4} 
+
4 \cdot \frac{29}{5760} \left< \; \cdot \; \right>_{H_{2,3}} 
+
4 \cdot \frac7{240} \left< \; \cdot \; \right>_{H_{2,2,2}}.
$$
The factors $8$, $4$, and $4$ are the numbers of symmetries
of the graphs.
\end{example}

Note that although from the point of view of the intersection
theory the case $\beta=1$ is the simplest and the most natural,
the graph representation of the corresponding bracket is
in no way simpler than for any other bracket. Note also
that the three simple graphs we have chosen can be
replaced by any other graphs with the same valencies
of vertices.

\begin{definition}
\label{Def:degree}
Let $\left< \; \cdot \; \right>$ be a bracket satisfying the
string and dilaton relations. 

The {\em genus} of the bracket is the number $g$ that appears in the
dilaton relations. 

We say that the bracket is {\em of pure degree $b$} if
it vanishes unless $b + \sum d_i = 3g-3+n = \dim \oM_{g,n}$. 
The bracket is {\em of nonnegative degree} if it is a finite linear
combination of brackets of pure nonnegative 
degrees, in other words, if it vanishes unless
$$
3g-3+n \geq  \sum d_i .
$$
\end{definition}

The bracket $\left< \beta \; \dots \; \right>$ is 
always of nonnegative degree and is of pure
degree $b$ if $\beta$ is a $2b$-cohomology class.

\begin{remark}
In some cases, Theorem~\ref{Thm:allbrackets} below can be 
generalized to brackets whose degree may be negative,
but bounded from below. To do that,
one must consider other types of marked graphs
then the ones we introduced in Definition~\ref{Def:decgraph}.
However, brackets with components of negative degree
do not naturally appear either in the study of
intersection theory of moduli spaces or in the
enumeration of ramified coverings. Therefore we
do not investigate them further.
\end{remark}

\begin{theorem} \label{Thm:allbrackets}
For any rational-valued bracket $\left< \; \cdot \; \right>$ 
of nonnegative degree, satisfying the string and
dilaton relations there exists a set of 
simple graphs $H_1, \dots, H_l$ and rational numbers
$q_1, \dots, q_l$ such that 
$$
\left< \; \cdot \; \right> \;
= \;
q_1 \, \left< \; \cdot \; \right>_{H_1} \; + \; \dots \; + \;
q_l \, \left< \; \cdot \; \right>_{H_l} \; .
$$
\end{theorem}

\paragraph{Proof.}
Let $\left< \; \cdot \; \right>$ be a genus $g$ bracket of
nonnegative degree. It follows from the inequality 
$3g-3+n - \sum d_i \geq 0$ that {\em if $d_i \geq 2$ for all $i$,
then $n \leq 3g-3$}. In other words, there is only a
finite number of brackets
that cannot be simplified using the string or the dilaton
relation. We will call the values of these brackets
the {\em initial values}. The initial values and the
string and dilaton relations completely determine all values
of the bracket.

Consider the set $S$ of all initial values of the bracket.
An element $v$ of $S$ is a monomial
$\tau_{d_1} \dots \tau_{d_m}$
(with $d_i \geq 2$ for all $i$) and a rational
number $q_v = \left< \tau_{d_1} \dots \tau_{d_m}\right>$.
For a given $v$, we denote by $b$ the number
$b = 3g-3+m - \sum d_i$ (the degree of the initial value).

To each initial value $v \in S$ we assign an arbitrary simple
graph $H_v$ with $m$ vertices of valencies $d_i+1$, $1 \leq i \leq m$, 
and a vertex $*$ of valency $b+g-1$. The condition
of nonnegative degree $b \geq 0$, together with $g \geq 2$ 
insures that the valency of the 
vertex $*$ is positive. It is easy to
see that a simple graph with these conditions can
always be constructed. Moreover, it can be chosen to be connected.
A simple computation shows that a simple graph like that
has $2g-2+m$ edges, in other words, it has Euler characteristic
$3-2g$.

We claim that the bracket $\left< \; \cdot \; \right>$
can be represented as
$$
\left< \; \cdot \; \right>
=
\sum_{v \in S} q_v\;  |\Aut(H_v)| \; \left< \; \cdot \; \right>_{H_v}.
$$
Indeed, the left-hand side bracket and the right-hand side
bracket have the same initial values and satisfy the
string and dilaton relations. Therefore they coincide.

Thus we have represented any given bracket 
$\left< \; \cdot \; \right>$ as a linear combination of
brackets assigned to simples graphs.
\qed

\begin{proposition} \label{Prop:graphs->A}
Consider a simple graph $H$. Let $e$
be its number of edges and $v$ its number of
vertices without counting the vertex $*$.  Denote by
$|\Aut(H)|$ the number of automorphisms of the
graph $H$. Then the series
$$
F_H(q) = 
\sum_{n \geq 0} \frac{q^n}{n!}
\sum_{d_1, \dots, d_n} 
\left< \tau_{d_1} \dots \tau_{d_n} \right>_H
$$
equals $\frac1{|\Aut(H)|} 
\, Y^v(1+Z)^e$ and, as a consequence,
lies in the algebra $\A$.
\end{proposition}

Various facts similar to this proposition were
independently discovered by D.~Zagier
(unpublished).

\begin{figure}[h]
\begin{center}
\
\epsfbox{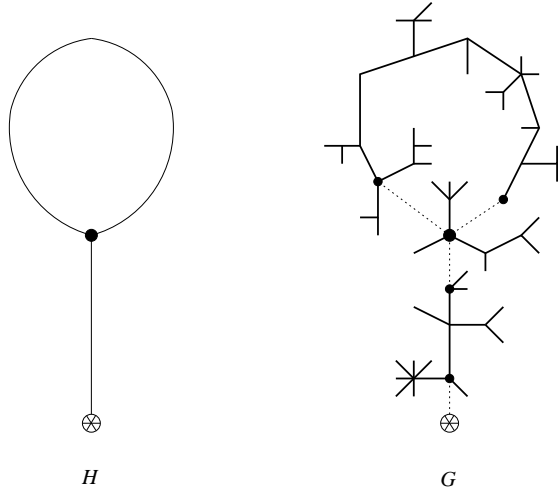}

\caption{Constructing the extensions of a simple graph $H$.}
\label{Fig:extension}
\end{center}
\end{figure}

\paragraph{Proof.} 
The series $F_H$ is the exponential generating series
that enumerates the $n$-decorated graphs $G$ whose
simplification equals $H$. Every graph $G$ like that
can be obtained by the following procedure shown in
Figure~\ref{Fig:extension}. First to every vertex of $H$
we assign a rooted tree. Second, to every
edge of $H$ we assign either a tree with two marked 
vertices (possibly twice the same vertex) or the empty tree
(with $0$ vertices).
Now join the marked vertices of these trees as shown
in the figure to obtain a graph $G$ homotopic to $H$.
Finally, number the vertices of $G$.

Since the exponential generating function for rooted
trees (respectively, for trees with two marked
vertices) is $Y$ (respectively, $Z$), we obtain
$$
F_H = \frac1{|\Aut(H)|} Y^v (1+Z)^e,
$$
where the factor $1/|\Aut(H)|$ comes from the fact
that we have counted as different the graphs $G$ obtained
from each other by symmetries of $H$.
Recalling that $Y$ and $Z$ are generators of $\A$, we see that
$F_H$ lies in $\A$.
This proves the proposition. \qed

Now we have all the necessary elements to prove 
Theorem~\ref{Thm:integral}, that we restate here.

\paragraph{Theorem 1}
{\em For any $g \geq 2$ and $\beta \in H^*(\oM_{g,0}, \Q)$, 
the power series
$$
F_{g,\beta} (q) =
\sum_{n \geq 0} \frac{q^n}{n!}
\int_{\oM_{g,n}} \frac{\beta}{(1-\psi_1) \dots (1-\psi_n)}
$$
lies in the algebra $\A$.}

\paragraph{Proof.} We have
$$
F_{g,\beta} (q) = \sum_n \frac{q^n}{n!}
\sum_{d_1 \dots d_n} \left< \beta \; \tau_{d_1} \dots \tau_{d_n} \right>
\, .
$$
According to Theorem~\ref{Thm:allbrackets}, the bracket
$\left< \beta \; \dots \; \right>$ is a linear combination
(with rational coefficients)
of brackets $\left< \; \cdot \; \right>_H$ for a finite
number of simple graphs $H$. The assertion of 
Theorem~\ref{Thm:integral} now follows from 
Proposition~\ref{Prop:graphs->A}. \qed

\subsection{A generalization}
\label{Ssec:generalization}

Here we formulate and prove a generalization of 
Theorem~\ref{Thm:integral}. The class $\beta$ will now
be a cohomology class of $H^*(\oM_{g,p}, \Q)$ for some
fixed nonnegative integer $p$. Moreover, the generalization
can be directly applied to the Ekedahl-Lando-Shapiro-Vainshtein
(ELSV) formula from~\cite{ELSV}, see  Section~\ref{Sec:coverings},
Theorem~\ref{Thm:ELSV}. 

Let $p \geq 0$ be an integer and $\beta \in H^*(\oM_{g,p}, \Q)$
a cohomology class of $\oM_{g,p}$. Fix $p$ positive integers
$b_1, \dots, b_p$ and denote by $r$ the number
$r = \sum (b_i-1)$.

\begin{theorem}
\label{Thm:gen}
For any $g,p$ such that $2-2g-p<0$ and for any positive
integers $b_1, \dots, b_p$, the power series
$$
F_{g;\beta;b_1, \dots, b_p} (q) =
$$
$$
\sum_{n \geq p+r} \frac{q^n}{(n-p-r)!}
\int\limits_{\oM_{g,{n-r}}} \frac{\beta}
{(1-b_1\psi_1) \dots (1-b_p\psi_p)
(1-\psi_{p+1}) \dots (1-\psi_{n-r})}
$$
lies in the algebra $\A$.
\end{theorem}

The proof goes along the lines of that of Theorem~\ref{Thm:integral},
but differs in some details. 

First, we slightly modify the
definition of an $n$-decorated graph.

\begin{definition}
An $n$-decorated graph $G$ is a graph with $n+1$ vertices
one of which is labeled by $*$ and the others are numbered
from $1$ to $n$. Moreover, we suppose that every connected
component of $G$ either contains one of the vertices
$*,1, \dots, p$ or has at least two independent cycles.
\end{definition}

\begin{definition}
Let $p$ be a nonnegative integer and $G$ an $n$-decorated graph. 
The $p$-{\em simplification of $G$} is the
graph $H$ obtained from $G$ by (i)~applying as many times as possible
the operations (S) and (D) without erasing the vertex $*$ or 
changing its valency, and also without erasing the vertices with numbers
from $1$ to $p$; (ii)~forgetting the numbers of the vertices,
except those from $1$ to $p$.

Conversely, the graph $G$ is called an {\em extension} of $H$.

A $p$-{\em simple} graph is
any finite graph with $p$ vertices numbered from $1$ to $p$,
one more vertex labeled by $*$, and possibly several
non-numbered vertices, each of which has a valency $\geq 3$. 
\end{definition}

It is easy to see that the $p$-simplification
of a graph $G$ is unique, in other words, it does not depend
on the order in which one applies the operations (S) and (D).

To any $p$-simple graph $H$ we assign a bracket
$\left< \; \cdot \; \right>_H$ using the same rule
as in Definition~\ref{Def:bracketH} replacing the
word ``simplification'' by ``$p$-simplification''.

\begin{proposition}
The bracket $\left< \; \cdot \; \right>_H$ for a
$p$-simple graph $H$ with Euler characteristic $3-2g$
satisfies the string and dilaton relations for $n \geq p$.
\end{proposition}

\paragraph{Proof.} Same as that of 
Proposition~\ref{Prop:bracketH}. \qed

\begin{proposition} \label{Prop:all}
Consider any rational-valued bracket either
of nonnegative degree and genus $g \geq 1$ or of
positive degree and genus $g=0$. Suppose that the
bracket satisfies the string and the dilaton relations for
$n \geq p$ ($p$ being a nonnegative integer). Then there
exists a set of $p$-simple graphs $H_1, \dots, H_l$
and of rational numbers $s_1, \dots, s_l$ such that
the bracket is equal to
$$
s_1 \left< \; \cdot \; \right>_{H_1} + \dots +
s_l \left< \; \cdot \; \right>_{H_l}.
$$
\end{proposition}

\paragraph{Proof.} The proof repeats that of 
Theorem~\ref{Thm:allbrackets}. The restriction
that the genus $g$ and the degree $b$ of the bracket
cannot vanish simultaneously comes from the
fact that the vertex $*$ must have the valency
$b+g-1$. \qed

\paragraph{Proof of Theorem~\ref{Thm:gen}.}
The generating series $F_{g,\beta;b_1, \dots, b_p}(q)$
from Theorem~\ref{Thm:gen} can be rewritten as
$$
F_{g,\beta;b_1, \dots, b_p}(q)
=
\sum_{n \geq p+r} \frac{q^n}{(n-p-r)!}
\sum_{d_1, \dots, d_n}
b_1^{d_1} \dots b_p^{d_p}
\left< \beta \; \tau_{d_1} \dots \tau_{d_{n-r}} \right> \; .
$$

The exceptional case $g=0$, $\deg \beta = 0$ can be
easily treated using the explicit values of the
bracket given in Section~\ref{Ssec:exceptions}.
We leave this as an exercise to the reader.

In all the other cases, by Proposition~\ref{Prop:all} the bracket 
$\left<\beta \; \dots \; \right>$ can be decomposed as
a linear combination over $\Q$ of brackets
$\left< \; \cdot \; \right>_H$ for some $p$-simple graphs
$H$. Therefore it is enough to prove that the series
$$
F_{g,H;b_1, \dots, b_p}(q)
=
\sum_{n \geq p+r} \frac{q^n}{(n-p-r)!}
\sum_{d_1, \dots, d_n}
b_1^{d_1} \dots b_p^{d_p}
\left< \tau_{d_1} \dots \tau_{d_{n-r}} \right>_H 
$$
lies in $\A$ for every $p$-simple graph $H$. This
series can be rewritten as follows
$$
F_{g,H;b_1, \dots, b_p}(q)
=
\sum_{n \geq p+r} \frac{q^n}{(n-p-r)!}
\sum_G b_1^{v_1-1} \dots b_p^{v_p-1} \; ,
$$
where the second sum is taken over all $(n-r)$-decorated 
graphs $G$ that simplify into $H$, while $v_1$, \dots
$v_p$ are the valencies of the first $p$ vertices
of $G$. We are going to give a way of enumerating
such graphs that shows that the corresponding
generating series lies in $\A$.

Consider a $p$-simple graph $H$. Denote
by $V_n$ the set $\{ 1, \dots, n \}$.
Consider the set $C_n$ of the following objects. 
(i)~An ordered list of
$p$ subsets $U_1, \dots, U_p \subset V_n$
satisfying $|U_i| = b_i$. (Here $|U_i|$ is the
number of elements in $U_i$.)
(ii)~A graph $\widehat G$ whose set of vertices
is $V_n \cup \{ * \}$. Moreover, we impose the
condition that if, in the graph $\widehat G$, 
we glue together the vertices of $U_i$ for each $i$
and attribute to the obtained vertex the label $i$,
we obtain an extension $G$ of the $p$-simple graph $H$.
Note that the vertices of $G$ are numbered incorrectly
(i.e., not form $1$ to $n-r$).

Denote by $|C_n|$ the number of objects like that.
We claim that 
$$
\frac{|C_n|}{n!} = \frac1{(b_1-1)! \dots (b_p-1)!} \; \cdot \;
 \frac1{(n-p-r)!}
\sum_G b_1^{v_1-1} \dots b_p^{v_p-1} \; ,
$$
where, as before, the sum is taken over the $(n-r)$-decorated
graphs that simplify into $H$. Indeed, the vertex number
$i$, $1 \leq i \leq p$, in the graph $G$ is separated
into $b_i$ different vertices in the graph $\widehat G$.
Thus every graph $G$ will be counted with weight
$b_1^{v_1} \dots b_p^{v_p}$. On the other hand,
each graph $\widehat G$ possesses 
$n!/b_1! \dots b_p!$ 
different numberings of the unmarked vertices, 
while the graph $G$ has only $(n-p-r)!$ of them.

\begin{figure}[h]
\begin{center}
\
\epsfbox{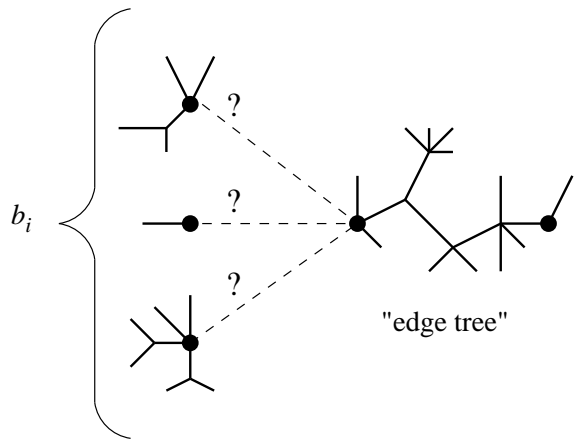}
\caption{There are $b_i$ choices for each edge of $H$
adjacent to the vertex number $i$.}
\label{Fig:choice}
\end{center}
\end{figure}

Now it remains to prove that the series
$$
\sum_{n \geq p+r} |C_n| \,  \frac{q^n}{n!}
$$
lies in $\A$. To do that, we propose the following way
to enumerate the objects from $C_n$. (i)~Construct
$b_1+ \dots + b_p$ disjoint rooted trees with vertices chosen
in the set $V_n$. Regroup these trees into $p$ forests
with $b_1$, \dots, $b_p$ trees. The roots of these trees
will form the sets $U_1$, \dots, $U_p$. (ii)~To each
non-numbered vertex of $H$ assign yet another rooted
tree with vertices in $V_n$. To each edge of $H$ assign
either an empty tree (with 0 vertices) or a tree
with two marked vertices. (iii)~Join the marked points
of the trees assigned to edges to the roots of the rooted trees
assigned to the non-numbered vertices (see
Figure~\ref{Fig:extension}).
(iv)~Now consider an edge of $H$ adjacent to a numbered
vertex, say to the vertex number $i$. There is a whole forest
of $b_i$ rooted trees assigned to this vertex. 
We choose one of these trees $T$ and join one of the marked points
of the ``edge tree'' to the root of $T$ 
(see Figure~\ref{Fig:choice}). Here it is crucial to note
that since the valency of each vertex of $H$ is finite
and fixed, {\em there is only a finite and fixed 
(i.e., independent of $n$) number
of choices}.

Our new description of the objects of $C_n$ shows
that the generating function 
$$
\sum_{n \geq p+r} |C_n| \,  \frac{q^n}{n!}
$$
is a finite sum (over a finite number of choices)
of finite products of series $Y$ (for rooted trees)
and $1+Z$ (for trees with two marked points).
Thus the generating function lies in $\A$. \qed

\section{Counting ramified coverings of the sphere}
\label{Sec:coverings}

This section is devoted to the enumeration of 
ramified coverings of the sphere by surfaces of a fixed genus $g$
and to a proof of Theorem~\ref{Thm:Hurwitz}.

\subsection{The ELSV formula}
\label{Ssec:ELSV}

Curiously, the most difficult part of the proof
of Theorem~\ref{Thm:Hurwitz} is the case with only
one complicated ramification point, $k=1$.
We know no other way to prove it
than to use the intersection theory
on moduli spaces and the results of the previous section.
Apart from these results, the main ingredient of the
proof is a theorem by 
T.~Ekedahl, S.~K.~Lando, M.~Shapiro, and
A.~Vainshtein that we formulate below after introducing
some notation.

Let $\mu= 1^{a_1} 2^{a_2} \dots$ be a partition
with degeneracy $r$. We define $|\Aut (\mu)|$ to be
$|\Aut (\mu)| = a_1! a_2! \dots$. For the formulation
of the theorem it is more convenient to switch to using
the additive notation for the partition $\mu$,
$\mu = (b_1, \dots, b_p)$, the $b_i$ being the
parts of $\mu$.

We will also need several cohomology classes on the
moduli space $\oM_{g,n-r}$. As before, $\cL_i$, $1 \leq i \leq n-r$
are the tautological line bundles (Notation~\ref{Not:MandL}),
and $\psi_i = c_1(\cL_i)$ are their first Chern classes. 
We also consider the {\em Hodge vector bundle} $W$. The fiber of
$W$ over a smooth curve is the set of holomorphic $1$-forms
on this curve. The fiber of $W$ over a general stable
curve is the set of global sections of its {\em dualizing sheaf}.
We do not give the details here (see the paper \cite{ELSV}
itself). Suffice it to note that $W$ is a vector bundle
of rank $g$. 

One property of the vector bundle $W$ will be important to us.
Recall that there is a canonical ``forgetful'' map from $\oM_{g,n-r}$
to $\oM_{g,0}$, forgetting the marked points and contracting the
components of the curve that have become unstable. The Hodge
bundle can be defined on both $\oM_{g,n-r}$ and $\oM_{g,0}$,
and {\em the former is the pull-back of the latter under the
forgetful map}. The same is true for the Chern classes
of the Hodge bundle, in particular for the
total Chern class $c(W^*)$ of its dual bundle.

The Hurwitz number $h_{g,n;\mu}$ is defined in 
Notation~\ref{Not:Hurwitz}. It counts the marked ramified
coverings with one ramification point of type $\mu$ and
$c(n) = 2n+2g-2-r$ simple ramifications.

Now we can write down the ELSV formula.

\begin{theorem} \label{Thm:ELSV} {\rm {\bf (The ELSV formula,}
\cite{ELSV}{\bf )}} \hspace{1em}
For any $g$, $n$, and $\mu$ such that $2-2g-(n-r) <0$, we have
$$
h_{g,n;\mu} \; = \;
\frac{(2n+2g-2-r)!}{|\Aut (\mu)|}
\prod_{i=1}^p \frac{b_i^{b_i}}{b_i!} \; \times
$$
$$
\times \;
\frac1{(n-p-r)!}
\int_{\oM_{g,n-r}} \frac{c(W^*)}
{(1-b_1 \psi_1) \dots (1-b_p\psi_p)
(1-\psi_{p+1}) \dots (1- \psi_{n-r})}
\; .
$$
\end{theorem}

Note that this formula seems to be waiting for someone to 
apply Theorem~\ref{Thm:gen}.

\subsection{Proof of Theorem~\ref{Thm:Hurwitz}}
\label{Ssec:ThmHurwitz}

Recall that $h_{g,n; \mu_1, \dots, \mu_k}$ is
the number of $n$-sheeted marked
ramified coverings of $\CP^1$ by a genus $g$ surface,
with $k$ ramification points of types $\mu_1, \dots, \mu_k$
and, in addition, $c(n) = 2n+2-2g-r$ simple
ramification points, where $r$ is the sum of degeneracies of
the partitions (see Notation~\ref{Not:Hurwitz}).
We are going to prove that all generating series for the numbers
$h_{g,n; \mu_1, \dots, \mu_k}$ with respect to the number
of sheets $n$ lie, once again, in the algebra $\A$.
In Section~\ref{Sec:conclusion} we explain why these
particular generating series are interesting to consider.

Let us restate the theorem we are going to prove.

\paragraph{Theorem 2}
{\em Fix any $g \geq 0$, $k \geq 0$. If $g=1$, we suppose that
$k \geq 1$. Then for any partitions
$\mu_1, \dots, \mu_k$, the series
$$
H_{g;\mu_1, \dots, \mu_k} (q) =
\sum_{n \geq 1} \frac{h_{g,n;\mu_1, \dots, \mu_k}}{c(n)!} 
\, q^n
$$
lies in the algebra $\A$.}

\paragraph{Proof.}
The theorem is proved by induction on the number $k$
of partitions.

{\bf Base of induction.} 
For $k=0,1$, the result is obtained by a direct application
of Theorem~\ref{Thm:gen} to the ELSV formula. 
In the case $k=0$, we must use
the ELSV formula with an empty partition $\mu$.

There are three exceptional cases in which Theorem~\ref{Thm:gen}
cannot be applied: $g=0$, $k=0$; $g=0$, $k=1$, $p \leq 2$; $g=1$,
$k=0$. These cases are discussed in Remark~\ref{Rem:exceptions}
below. It turns out that the assertion of Theorem~\ref{Thm:Hurwitz}
fails only if $g=1$, $k=0$, as stated in the formulation.

{\bf Step of induction.}
The step of induction is an almost exact repetition of 
the proof of Theorem~2
from our previous work~\cite{Zvonkine}. We only give a short summary of
the argument here, referring to~\cite{Zvonkine} for details.

The proof goes in the spirit of~\cite{GouJac}. 
Recently, M.~Kazaryan introduced
a semi-group of colored permutations, which seems to be best
fit for describing the proof (work in preparation).

It is easy to see that there is only a finite number of 
possible cycle structures for a permutation that can be
obtained as a product of two permutations with given cycle
structures $\mu_1$ and $\mu_2$. 

Let $\mu_1$ and $\mu_2$ be two partitions from the list
$\mu_1, \dots, \mu_k$. We can move the two corresponding 
ramification points on $\CP^1$ towards each other until
they collapse. We obtain a new (not necessarily connected)
ramified covering. Its monodromy at the new ramification point
is the product of the monodromies of the two points that have
collapsed. 

Let us choose one
of the possible cycle structures of the product monodromy
and also one of the possible ways in which the 
covering can split into connected components.
By the induction assumption, we obtain a series from
the algebra $\A$ assigned to each connected component of
the covering. Indeed, each connected component is itself
a ramified covering of the sphere as in Theorem~\ref{Thm:Hurwitz},
but with $k-1$ fixed ramification types instead of $k$.
We obtain the generating series for the number of nonconnected 
ramified coverings by multiplying the series that correspond
to the connected components. Since it is a finite product
of series lying in $\A$, we obtain again a series from $\A$.

Finally, we must add the generating series described above
for all choices of types of nonconnected coverings. Since the
number of choices is finite, we obtain, once again,
a series from $\A$.
\qed

\begin{remark} \label{Rem:exceptions}
Let us consider the exceptional cases $g=0$, $k=0,1$
and $g=1$, $k=0$. 

In the genus zero case, the ELSV formula transforms into
a much simpler Hurwitz formula~\cite{Hurwitz,Strehl},
which turns out to be applicable even if the multiple
ramification point has only $1$ or $2$ preimages.

We have, using the notation of Theorem~\ref{Thm:ELSV} and
Notation~\ref{Not:Hurwitz},
$$
h_{0,n;\mu} =
\frac{(2n-2-r)!}{|\Aut(\mu)|} \;
\prod_{i=1}^p \frac{b_i^{b_i}}{b_i!} 
\; \cdot \; 
\frac{n^{n-r-3}}{(n-p-r)!}.
$$
This formula is true for any $n \geq p+r$ and for
any partition $\mu$ (including even the empty partition).
We see that the corresponding generating series always
lies in the algebra $\A$.

The case $g=1$, $k=0$ is covered by the ELSV formula
with an empty partition $\mu$. Consider the moduli 
space $\oM_{1,1}$. As in the last paragraph of 
Section~\ref{Ssec:exceptions}, denote by $\beta$ the $2$-cohomology
class of $\oM_{1,1}$ whose integral over the fundamental
homology class equals~1. One can prove that the Hodge
bundle over $\oM_{1,1}$ is a line bundle with first
Chern class $\beta/24$. Therefore we obtain
$$
h_{1,n;\emptyset} = (2n)! \cdot
\frac{1}{n!} \int_{\oM_{1,n}}
\frac{1-\frac{1}{24}\beta}{(1-\psi_1) \dots (1-\psi_n)} \; .
$$
Using the calculations of Section~\ref{Ssec:exceptions},
we obtain
$$
\sum_{n \geq 1} \frac{h_{1,n;\emptyset}}{(2n)!} \, q^n
= 
\frac{1}{24} \, \sum_{n \geq 1} \frac{A_n}{n} \frac{q^n}{n!} \; .
$$
This series does not lie in $\A$ (and constitutes the
only exception to the general rule). It suffices
to consider the partition $\mu = (1)$, which
amounts to distinguishing one sheet in the ramified
covering, to obtain the series
$$
\frac{1}{24} \, \sum_{n \geq 1} 
\frac{A_n}{n!} \, q^n \; \in \; \A \, .
$$
\end{remark}

\section[Random metrics versus random 
differentials]{Random Riemannian metrics versus random 
abelian differentials}
\label{Sec:conclusion}

In this section we do not prove any theorems, but
discuss the relation between the enumeration of
ramified coverings of the sphere and the 2-dimensional
gravity. We also draw a parallel between the
study of spaces of Riemannian metrics using
ramified coverings of the sphere and the study of
spaces of abelian differentials using ramified coverings
of the torus.

\subsection{Two models of 2-dimensional gravity}

In every problem of statistical physics one starts
with introducing a space of states and by assigning
an energy to every state. 

In 2-dimensional gravity, a {\em state} is a
2-dimensional compact oriented real
not necessarily connected surface endowed with a
Riemannian metric. Two surfaces like that are
equivalent, i.e., correspond to the same state,
if they are isometric.

Consider a surface $S$ with a Riemannian metric. Let
$\chi(S)$ be its Euler characteristic and $A$ its
total area. To such a surface one assigns an {\em energy}
$$
E = \lambda A + \mu \chi(S).
$$
Here $\lambda$ and $\mu$ are two constants called
the cosmological constant and the gravitational
constant, respectively. Note that $\chi(S)$ is
actually the integral over $S$ of the curvature
of the metric. The fact that this integral takes
such a simple form is special to dimension~2.

Now the first thing to do is to compute the
{\em partition function} 
$$
Z(\lambda, \mu) = \int_{\mbox{states}} e^{-E},
$$
or, equivalently, the {\em free energy}
$$
F(\lambda, \mu) = \ln Z(\lambda, \mu)
= \sum_{g \geq 0} \int_{\mbox{metrics}}
e^{-E}.
$$
The free energy is the sum of contributions of
connected surfaces, while the partition function
is the sum of contributions of all surfaces.

Neither of the above integrals is well-defined
mathematically, but we would still like to
compute them. To do that, physicists introduced
a {\em discrete model} of Riemannian metrics,
replacing them by quadrangulations~\cite{BreKaz,Witten}
(see also~\cite{thebook}, Chapter~3 
for a mathematical description). In this model,
instead of considering Riemannian metrics, one considers
metrics obtained by gluings of squares of area $\varepsilon$.
Our goal is to convince the reader that the ramified coverings
of the sphere provide a new (maybe more natural)
discrete model of Riemannian metrics.

Fix a positive number $\varepsilon$. Consider a sphere
with the standard (round) Riemannian metric
and with total area $\varepsilon$. On this sphere,
choose at random $2n+2g-2$ points. Now chose
a random connected $n$-sheeted covering of the sphere with simple
ramifications over the $2n+2g-2$ chosen points. 
The covering surface $S$ will automatically be of genus $g$. 
The metric on the sphere can be lifted to $S$, which will
give us a metric with constant positive curvature except
at the critical points, where it has conical singularities
with angles $4\pi$. This metric is, of course, not Riemannian.
However, one can argue that if $\varepsilon$ is very small and the number
of sheets very large, a random metric obtained
in this way looks similar to a random
Riemannian metric (unless we look at them through a
microscope to reveal the difference). We do not know
any rigorous statement that would formalize this
intuitive explanation, but the same argument is
used by physicists to justify the usage of quadrangulations.

Using our discrete model of metrics, one can
write the free energy for the 2-dimensional
gravity in the following way:
$$
F(\lambda, \mu) = 
\sum_{g,n} \frac{\varepsilon^{2n+2g-2}}{(2n+2g-2)!} \; 
h_{g,n; \emptyset} \; e^{-n\lambda \varepsilon - \mu (2-2g)}.
$$
Here $\varepsilon^{2n+2g-2}/(2n+2g-2)!$ is the volume of the space
of choices of $2n+2g-2$ unordered points on the sphere
of area $\varepsilon$.

Using Theorem~\ref{Thm:Hurwitz} (and some additional
considerations to find the degree in $Y$ and $Z$ of
the series that are involved), one can show that
the leading term of the asymptotic of 
$h_{g,n; \emptyset}/(2n+2g-2)!$ 
as $n$ tends to infinity (while $g$ is fixed) looks like
$$
\frac{h_{g,n; \emptyset}}{(2n+2g-2)!} 
\sim e^n n^{\frac52 (g-1) - 1} b_g,
$$
with some constants $b_g$. The first values of these constants are
$$
b_0 =\frac1{\sqrt{2\pi}}, \quad
b_1 =\frac{1}{2^4 \cdot 3}, \quad
b_2 = \frac1{\sqrt{2\pi}} \; \frac{7}{2^5 \cdot 3^3 \cdot 5},
$$
$$
b_3 = \frac{5 \cdot 7^2}{2^{16} \cdot 3^5}, \quad
b_4= \frac1{\sqrt{2\pi}} \;
\frac{7 \cdot 5297}{2^{11} \cdot 3^8 \cdot 5^2 \cdot 11 \cdot 13}.
$$

Now we make the final step by letting $\varepsilon$ tend
to $0$ in the expression for the free energy $F$. 
To obtain an interesting limit for the free energy, 
we must {\em make $\lambda$ and $\mu$ depend on $\varepsilon$}.
We want to use
$$
\sum_{n \geq 1} n^{\alpha-1} e^{-\varepsilon n}  \sim 
\frac1{\varepsilon^\alpha} \Gamma(\alpha)
\quad \mbox{ as } \varepsilon \rightarrow 0.
$$
Therefore we let
$\lambda \varepsilon -2 \ln \varepsilon - 1 \rightarrow 0$, 
while
$$
y = 
\frac{(\varepsilon e^\mu)^{4/5}}
{\lambda \varepsilon -2 \ln \varepsilon-1}
$$
remains fixed. This gives us the final expression
of the free energy, now depending on only one variable
$y$:
$$
F(y) = \Gamma (-5/2) \, b_0 \, y^{5/2} - b_1 \ln y
+ \sum_{g \geq 2} 
\Gamma \biggl(5(g-1)/2\biggr) \,
\, b_g \cdot y^{5(1-g)/2}.
$$
The coefficients
$\Gamma (5(g-1)/2) \, b_g $
are rational for odd $g$ and rational multiples
of $\sqrt2$ for even $g$.

\bigskip

Our above treatment is parallel to the treatment of
the quadrangulation model in~\cite{Witten}. 
Denote by $Q_{g,n}$ the number of ways to
divide a surface of genus $g$ into $n$ squares.
Then the study of the quadrangulation model
involves the asymptotic of $Q_{g,n}$, which is
given by
$$
Q_{g,n} \sim 12^n n^{\frac52(g-1) -1} b'_g,
$$
for another sequence of constants $b'_g$. This
sequence was studied using matrix integrals, and it
is known that a generating function for the sequence
$b'_g$ satisfies the Painlev\'e I equation. Our 
numerical experiments show that the expressions for the
free energy obtained in both models coincide up to
a rescaling of $y$. More precisely, we formulate
the following conjecture.

\begin{conjecture}
We have
$$
b'_g = 2^{\frac32 (g-1)+1} \cdot b_g .
$$
An equivalent statement: 
the function $u(y) = F''(y)$ satisfies the
Painlev\'e I equation
$$
\frac16 u''(y) + u(y)^2 = 2y.
$$
\end{conjecture}

It is also an open problem to find the Korteweg--de~Vries
hierarchy using the enumeration of ramified coverings,
as it was done with quadrangulations in~\cite{Witten}.

We believe that the connection between the two models
can be obtained using A.~Okounkov's results on random
Young diagrams: the identity between the distribution of highest
eigenvalues of a random hermitian matrix and the distribution
of lengths of longest columns in a random Young 
diagram~\cite{Okounkov1}; and an appearance of
integrable hierarchies in the study of random Young
diagrams~\cite{Okounkov2}.

\subsection{Ramified coverings of a torus and abelian
differentials}

Fix an integer $g \geq 1$ and a list of $p$ nonnegative
integers $b_1, \dots, b_p$ with the condition $\sum b_i = 2g-2$. 
We consider the space
$D_{g;b_1, \dots, b_p}$ of abelian differentials
on Riemann surfaces of genus $g$, with zeroes of
multiplicities $b_1, \dots, b_p$. More precisely,
$D_{g;b_1, \dots, b_p}$ is the space of triples
$(C, \{ x_1, \dots, x_p \} , \alpha)$, where
$C$ is a smooth complex curve, $x_1, \dots, x_p \in C$
are distinct marked points, and $\alpha$ is an
abelian ($=$ holomorphic) differential on $C$ whose
zero divisor is precisely $b_1 x_1 + \dots + b_p x_p$.

It turns out that the space $D_{g;b_1, \dots, b_p}$
has a natural {\em integer affine structure}.
This means that it can be covered by charts of
local coordinates in such a way that the change of
coordinates, as one goes from one chart to another,
is an affine map with integer coefficients. Such 
local coordinates are introduced in the following
way. Fix a basis $l_1, \dots, l_{2g+p-1}$ of 
the relative homology group
$H_1(C,\{x_1, \dots, x_p \}, {\mathbb Z})$.
Then the integrals of $\alpha$ over the cycles
$l_i$ are the local coordinates we need.
The {\em area function}
$$
A: (C, \alpha) \mapsto \frac{i}{2}
\int_C \alpha \wedge {\bar \alpha}
$$
is a quadratic form with respect to the affine
structure.

The integer affine structure allows one to define
a volume measure on the space
$D_{g;b_1, \dots, b_p}$. It is then a natural
question to find the total volume of the part
of the space $D_{g;b_1, \dots, b_p}$ defined
by $A \leq 1$ (the volume of the whole space
being infinite).

A.~Eskin and A.~Okounkov~\cite{EskOko} obtained an
effective way to calculate these volumes using
the asymptotic for the number of ramified coverings
of a torus. Consider the elliptic curve obtained
by gluing the opposite sides of the square
$(0,1,i,1+i)$ endowed with the abelian differential
$dz$. Given a ramified covering of this elliptic curve
with critical points of multiplicities
$b_1, \dots, b_p$, we can lift the abelian differential
to the covering curve and obtain a point of
$D_{g;b_1, \dots, b_p}$. One can then easily show that 
such points are densely and uniformly distributed in
$D_{g;b_1, \dots, b_p}$ if one considers coverings
with a big number of sheets. Moreover, 
R.~Dijkgraaf~\cite{Dijkgraaf} and S.~Bloch and
A.~Okounkov~\cite{BloOko} showed that the generating series for
the ramified coverings of the torus that arise
in this study are quasi-modular forms. In other
words, they lie in the algebra
$$
\Q [ E_2, E_4, E_6],
$$
where $E_{2k}$ are the Eisenstein series
$$
E_{2k}(q) = \frac12 \zeta(1-2k) + \sum_{n \geq 1}
\left( \sum_{d | n} d^{2k-1} \right) q^n.
$$

We conclude with the following comparison between the
counting of ramified coverings of a sphere and of a torus.

\bigskip 

Sphere: The generating series enumerating the ramified
coverings lie in the algebra $\A$.

Torus: The generating series enumerating the ramified
coverings lie in the algebra of quasi-modular forms.

\bigskip

Sphere: The coefficients of a
generating series grow as $e^n \cdot n^{\gamma} \cdot c$.
The exponent $\gamma$ is a half-integer. The
number $1- \gamma$ is called the
{\em string susceptiblity}.
The constant $c$ is an observable in 2-dimensional gravity.

Torus: The  sum of the first $n$ coefficients of a
generating series grows as $n^d \cdot c$. The number
$d$ is the complex dimension of the corresponding space of
abelian differentials. The constant $c$ is its volume.

\bigskip

Sphere: It is known that some particular observables
in 2-dimensional gravity can be arranged into a generating
series (in an infinite number of variables), that turns out
to be a $\tau$-function for the Korteweg--de~Vries hierarchy.
These observables have not been found yet in the model
with coverings of the sphere.

Torus: As far as we know, nobody has tried to arrange the
volumes of the spaces of abelian differentials into
a unique generating series.

\end{document}